\documentstyle[12pt]{article}

\input amssym.def
\input amssym.tex

\ifx\shlhetal\undefinedcontrolsequence\let\shlhetal\relax\fi
\font\teneuf=eufm10
\font\seveneuf=eufm7
\font\fiveeuf=eufm5
\newfam\frakturfam
\textfont\frakturfam\teneuf
\scriptfont\frakturfam\seveneuf
\scriptscriptfont\frakturfam\fiveeuf

\newcommand{\bbn}{{\Bbb N}}

\newcommand{\Prob}{{\rm Prob}}
\newcommand{\Lim}{{\lim}}

\newcommand{\dom}{{\rm dom}}
\newcommand{\lesdot}{\mathrel{\mathord{<}\!\!\raise 0.8
pt\hbox{$\scriptstyle\circ$}}}   

\renewcommand{\ln}{{\rm ln}}
\newcommand{\cY}{{\cal Y}}
\newcommand{\cM}{{\cal M}}
\newcommand{\bu}{{\bf u}}
\newcommand{\bc}{{\bf c}}

\newtheorem{theorem}{Theorem}[section]

\newtheorem{discussion}[theorem]{Discussion}

\title{ Very weak zero one law for random graphs with order and random binary
functions }
\author{{\bf Saharon Shelah}\thanks{\ \ The research partially
supported by the United States-Israel Binational Science Foundation
and  NSF under grant \#144-EF67;\quad Publication no 548.}\\
Institute of Mathematics, The Hebrew University\\
Department of Mathematics, Rutgers University\\
Department of Mathematics, University of WI, Madison
}
\date{done: August 1993, October 31, 1993\\
last corrections introduced June 5, 1996\\
printed: \today}

\begin{document}
\maketitle
\vfill
\eject

\section{Introduction}
Let $G_<(n,p)$ denote the usual random graph $G(n,p)$ on a totally
ordered set of $n$ vertices.  (We naturally think of the vertex
set as $1,\ldots,n$ with the usual $<$).  We will fix $p=\frac{1}{2}$
for definiteness.
Let $L^<$ denote the first order
language with predicates equality $(x=y)$, adjacency $(x\sim y)$
and less than $(x<y)$.  For any sentence $A$ in $L^<$ let
$f(n)=f_A(n)$ denote the probability that the random $G_<(n,p)$
has property $A$.  It is known Compton, Henson and Shelah
\cite{CHSh245} that there are $A$ 
for which $f(n)$ does not converge.  Here we show what is called
a {\em very weak zero-one law} (from [Sh 463]):
\begin{theorem}\label{Thm1}  For every $A$ in language $L^<$
\[
\lim_{n\rightarrow\infty} (f_A(n+1)-f_A(n)) = 0 \]
\end{theorem}
Note, as an extreme example, that this implies the nonexistence
of a sentence $A$ holding with probability $1-o(1)$ when $n$ is
even and with probability $o(1)$ when $n$ is odd (as in Kaufman,
Shelah \cite{KfSh201}).
\par In \S 2 we give the proof, based on a circuit complexity result.
In \S 3 we prove that result, which is very close to the now classic
theorem that parity cannot be given by an $AC^0$ circuit.  In \S 4
we give a very weak zero-one law for random two-place
functions.  The proof
is very similar, the random function theorem being perhaps of
more interest to logicians, the random graph theorem to
discrete mathematicians.

The reader should thank Joel Spencer who totally rewrote the paper
(using the computer science jargon rather than the logicians one), and
with some revisions up to the restatement in the proof of
2.1 but with 3.1, this is the version presented here. We
thank the referee for comments on the exposition, and we thank Tomasz
\L uczak and Joel Spencer for reminding me this problem on $G_{<}(n,
p)$ in summer 93.

On a work continuing this of Boppana and Spencer see \ref{disc1.6}(5).

\section{The Proof.}
Let $G$ be a fixed graph on the ordered set $1,\ldots,2n+1$.  For a
property $A$ and for $i=n,n+1$ let $g(i)=g_{G,A}(i)$ denote the
probability that $G\restriction_S$ satisfies $A$ where $S$ is chosen uniformly
from all subsets of $1,\ldots,2n+1$ of size precisely $i$.  We shall
show
\begin{theorem}\label{Thm2}  $g(n+1)-g(n)=o(1)$ \end{theorem}
More  precisely, given $A$ and $\epsilon>0$ there exists $n_0$ so that
for any $G$ as above with $n\geq n_0$ we have $|g(n+1)-g(n)|<\epsilon$.
\par We first show that Thm.\ref{Thm1} follows from Thm.\ref{Thm2}.  The idea
is that a random $G_<(i,p)$ on $i=n$ or $n+1$ vertices is created by
first taking a random $G_<(2n+1,p)$ and then restricting to a random
set $S$ of size $i$.  Thus (fixing $A$) $f(n),f(n+1)$ are the averages 
of $g_G(n),g_G(n+1)$ over all $G$.  By Thm.\ref{Thm2} we have
$g_{G,A}(n)-g_{G,A}(n+1)=o(1)$ 
for all $G$ and therefore their averages are only $o(1)$ apart.
\par Now we show Thm.\ref{Thm2}.  Fix $G,A$ as above.  Let $P(S)$ be the
Boolean value of the statement that $G|_S$ satisfies $A$.
For $1\leq x\leq 2n+1$
let $z_x$ denote the Boolean value of ``$x\in S$''  so that $P(S)$ is
a Boolean function of $z_1,\ldots,z_{2n+1}$.  We claim this function
has a particularly simple form.  Any $A$ can be built up from
primitives $x=y,x<y,x\sim y$ by $\wedge,\neg$ and, critically,
$\exists_x$.  As $G$ is fixed the primitives have values true
or false.  Let $\wedge,\neg$ be themselves.  Consider $\exists_x W(x)$
where for each $1\leq x\leq 2n+1$ we let $W(x)$ on $G|_S$ is given by $W^*(x)$.
Then $\exists_x W(x)$ has the interpretation $\exists_{x\in S} W(x)$
which is expressed as $\vee_{x=1}^{2n+1} (z_x\wedge W^*(x))$.  For
convenience we can be redundant and replace $\forall_x W(x)$ by
$\wedge_{x=1}^{2n+1} (z_x\Rightarrow W^*(x))$.  For
example $\forall_x\exists_y x\sim y$ becomes
\[ \wedge_x \left[ z_x \Rightarrow \vee_{y\sim x} z_x\right]  \]
Thus $P(S)$ can be built up from $z_1,\ldots,z_{2n+1}$ be means of 
the standard $\neg,\wedge,\vee$ and $\wedge,\vee$ over (at most)
$2n+1$ inputs.  That is (see \S 3) $P(S)$ can
be expressed by  an $AC^0$ circuit over $z_1,\ldots,z_{2n+1}$ (of
course with the number of levels bounded by the length $d_A$ of the
sequence $A$ (can get less) and the number of nodes bounded by $d_A n^{d_A}$).
Now $g(i)$, for $i=n,n+1$, is the probability $P$ holds when a
randomly chosen set of precisely $i$ of the $z$'s are set
to True.  From Thm.\ref{Thm3} below $g(n+1)-g(n)=o(1)$ giving
Thm. \ref{Thm2} and hence Thm. \ref{Thm1}.
\section{$AC^0$ Functions}
We consider Boolean functions of
$z_1,\ldots,z_m$.  (In our application $m=2n+1$.)  The functions
$z_i,\neg z_i$, called literals,  are the  level $0$ functions.  A level $i+1$
function is the $\wedge$ or $\vee$ of  polynomially many
level $i$ functions.  An $AC^0$ function is a level $d$ function
for any constant $d$.  By standard technical means we can express
any $AC^0$ function in a ``levelled'' form so that the level
$i+1$ functions used are either all $\wedge$s of level $i$ functions 
or all $\vee$s of level $i$ functions and the choice alternates with
$i$ (at most doubling the number of levels).
It is a classic result of circuit 
complexity that parity is not an $AC^0$ function.  Let $C$ be an
$AC^0$ function.
For $0\leq i\leq m$ let $f(i)=f_C(i)$ denote
the probability $C$ holds when  precisely $i$ of the $z_j$ are set to
True and these $i$ are chosen randomly.
\begin{theorem}\label{Thm3}
 $f(n+1)-f(n)=o(1)$
\end{theorem}

Called a restriction $\rho$ balanced if $|\{ i: \rho(i)=0\}|= |\{i:
\rho(i)=1\}|$.

Now more fully the theorem says
\begin{description}
\item[$(\ast)$] for every $\varepsilon$, $d$, $t$ there is
$n_{\varepsilon, d, t}$ satisfying: if $n\geq n_{\varepsilon, d, t}$ and
$C$ is an $AC^0$ Boolean circuit of $z_1, \ldots, z_{2n+1}$ of level
$\leq d$ with $\leq n^t$ nodes then $|f_C(n+1) - f_C(n)|< \varepsilon$.
\end{description}
This statement is proved by induction on $d$.

We choose the following
\begin{description}
\item[(i)] $c_0=(\ell n4)t> 0$
\item[(ii)] $\varepsilon=\frac{1}{2}$,
$\varepsilon_\ell=\frac{1}{2^{1+\ell}}$
\item[(iii)] $k$ is such that $\varepsilon \cdot k\geq t$
\item[(iv)] we choose $k_\ell$ inductively on $\ell\leq k$ such that $k_\ell$
large enough.
\item[(v)] $c_1$ a large enough real
\item[(vi)] $n_0$ is large enough
\end{description}
For a node $x$ of the circuit $C$ let ${\cal Y}_x$ be the set of nodes
which fans into it; (without loss of generality in the
level 1 we
have only OR).

First we assume $d> 2$.
Note
\begin{description}
\item[$\otimes_1$] drawing as below a balanced restriction $\rho$ with
domain with $\leq n$ elements, with probability $\geq 1-\varepsilon/3$ we
have: in $C^1=C\restriction_\rho$, every node of the level 1 (i.e.
for which $\cY_x$ is a set of atoms) satisfies $|\cY_x|\leq c_0 (\ln
n)$.
\end{description}
[Why? Choose randomly a set $\bu_0$ of $[n/2]$ pairwise disjoint pairs
of numbers among $\{1, \ldots, 2n+1\}$, and then for each $\{i, j\}\in
\bu$ decide with probability half that $\rho(i)=0$, $\rho(j)=1$ and with
probability half that $\rho(i)=1$, $\rho(j)=0$ (independently for
disjoint pairs). This certainly gives a balanced $\rho$.

Now if $x$ is a node of $C$ of the  level 1, the probability that
$\rho$ does not decide the truth value which the node compute is $\leq
(\frac{1}{4})^{|\cY_x|}$. Note: after drawing $\bu$, if $\cY_x$ contains
a pair from $\bu$ the probability is zero, we only increase compared to
drawing just a restriction. So the probability that for {\em
some} $x$ of the level 1 of $C$, $|\cY_x|\geq (\ln 4)t(\ln n)+1$ and
the truth value is not computed, is $\leq
|C|\times\big({1\over 4}\big)^{(\ell n 4)t(\ell n n)+1}
\leq 1/2$, so there is $\rho_0$ for
which for any such $x$ the truth value is computed.]

Next, we say that a restriction $\rho'$ extends a restriction $\rho$ if

\noindent $\rho'(i)\neq \rho(i) \Rightarrow \rho(i)=\ast$.
Now

\begin{description}
\item[$\otimes_2$] Choosing randomly a restriction $\rho_1$ as below we
have: $\rho_1$ is a balanced restriction extending $\rho$ such that
$|\{i: i\in \{1, \ldots, 2n+1\}, \rho(i)=\ast\}|\geq 2[n^\varepsilon]+1$
and with probability $\geq 1-\varepsilon/3$ for every node $y$ of $C$
of the  level 1 we have, $|\cY_y|\leq k$.
\end{description}
[Why? We draw a set $\bu_1$ of
$(2n+1-|\dom(\rho_0)|-(2[n^{\varepsilon_0}]+1))/2$ pairs from $\{i:
\rho_0(i)=\ast\}$ pairwise disjoint and for each $\{i, j\}\in \bu$,
decide with probability $\frac{1}{2}$ that $\rho_1(i)=0$, $\rho_1(j)=1$
and with probability half that $\rho_1(i)=1$, $\rho_1(j)=0$.

For each node $y\in C^1$ of the level 1 the probability that
``the
number of $y'\in \cY_y$ not assigned a truth value by
$\rho_1$ is $\geq k+1$'' is at
most $\left( \begin{array}{l} |\cY_y|\\ k+1 \end{array}\right) \times
\left(\frac{1}{2n^{\varepsilon_0}+1}\right)^{k+1}\leq (c_0\ln n)^{k+1}\cdot
n^{-\varepsilon_0(k+1)}< n^{-t}$.]

We now choose by induction on $\ell\leq k$ a restriction $\rho_{2,
\ell}$ such that
\begin{description}
\item[$\otimes_3$\ (a)] $\rho_{2, \ell_0}=\rho_1$, $\rho_{2,
\ell}\subseteq \rho_{2, \ell+1}$, $2n+1 - (2[n^{\varepsilon_\ell}]+1) =
|\dom \rho_{2, \ell}|$
\item[\quad\ \ (b)] every $y\in C$ of the  level 2 there is a set
$w_{y, \ell}$ of $\leq k_\ell$ atoms such that:
if $z\in \cY_y$, then $|\cY_z\setminus w_{y, \ell}|\leq k-\ell$.
\end{description}

Now for $C\restriction\rho_{2,k}$ we can invert AND and OR (multiplying
the size by a constant $\leq c_1$) decreasing $d$ by one
thus carrying the induction step.

For $\ell=0$ let $\rho_{2, 0}=\rho_1$. For $\ell+1$, for each $y\in C$
of level 2 let $\Xi= \{\nu: \nu\mbox{ a restriction with domain
}w_{y, \ell}\}$ let
$$
\cY^\nu_y=\{z\in \cY_y: \mbox{ the
truth value at $z$ is still not computed
under }\rho_{2, \ell}\cup \nu\},
$$
and try to choose by
induction on $i$ an atom
$z_{y, \ell, \nu, i}\in \cY^\nu_y\setminus \{z_{y, \ell, \nu, j}: j< i\}$,
such that $\dom(z_{y, \ell, \nu, i})$ is disjoint to $\bigcup\limits_{j<
i}\dom(z_{y, \ell, \nu, j})\setminus w_{y, \ell}$.
Let it be defined if $i< i_{y, \ell}$.

Now $\rho_{2, \ell+1}$ will for each $\nu\in\Xi$
decide that $\nu$ make the truth value
computed in $y$ true, or will leave only $\leq (k_{\ell+1}-k_\ell)/
2^{k_\ell}$ of the atoms in $\bigcup\limits_i \dom z_{y,
\ell, \nu,
i}\setminus w_{y, \ell}$ undetermined (this is done
as in the previous two stages).

But now by $\otimes_1$ + $\otimes_2$,
$C\restriction_{\rho_{2,k}}$ can be
considered having $d-1$ levels (because, as said above we
can invert the AND and OR in level 1 and 2).

We have translate our problem to one with
$[n^{\varepsilon_k}]$, $d-1$,
$\varepsilon_k(t+\varepsilon_1)$, $\frac{\varepsilon}{3}$ instead $n$,
$d$, $t$, $\varepsilon$ (the $t+\varepsilon$ is just for
$n^{t+\varepsilon}>c_1n^t$).

Also note: $\varepsilon$, $c_1$ does not depend on $n$. So we can use
the induction hypothesis. We still have to check the case $d\leq 2$, we
still are assuming  level 1 consist of cases of OR, and for
almost all random $\rho_1$ (as in $\otimes_1$) for every $x$ of
 level 1 we have $|\cY_x|\leq \bc_0 \ln n$ (so again changing $n$).

So as above we can add this assumption. Choose randomly a complete
restriction $\rho^0$ with $|\{i: \rho^0(i)=1\}|=n$, and let $\rho^1$ be
gotten from $\rho^0$ by changing one zero to $1$, so $|\{i:
\rho^1(i)=1\}|=n+1$.

Now the probability that $C\restriction_{\rho^0}=0$ but
$C\restriction_{\rho^1}=1$ is small: it require that for some node $x$
of level $1$ is made false in $C\restriction_{\rho^0}$ while there is no
such $x$ for $C\restriction_{\rho^1}$, but if $x(*)$ is such for
$C\restriction_{\rho^0}$ it is made true then with probability $\geq
1-\frac{|\cY_x|}{2n+1}\geq 1-\frac{\bc_0 \ln n}{n}$ the $z_i$ changed is
not in $\cY_{x(*)}$. Contradiction, thus finishing the
proof.


\section{Two Place Functions}
Here we consider the random structure $([n],F_n)$ where
$F_n(x,y)$ is a random function from $[n]\times[n]$ to
$[n]$.  (We no longer have order.  A typical sentence would
be $\forall_x\exists_y F(x,y)=x$.):  Again for any sentence $A$
we define $f(n)=f_A(n)$ to be the probability $A$ holds in
the space of structures on $[n]$ with uniform distribution.
Again it is known \cite{CHSh245} that convergence fails, there
are $A$ for which $f(n)$ does not converge.  Again our result
is a very weak zero one law.
\begin{theorem}\label{Thm4}
For every $A$
\[ \lim_{n\rightarrow\infty} f_A(n+1)-f_A(n) = 0 \]
\end{theorem}
Again let $m=2n+1$.  Let $F^*(x,y,z)$ be a {\em three}-place
function from $[m]\times[m]\times[m]$ to $[m]$.  For 
$S\subset [m]$ of cardinality $i=n$ or $n+1$ we define $F^*_S$,
a partial function from $[S]\times[S]$ to $[S]$ by setting
$F^*_S(x,y)=F^*(x,y,z)$ where $z$ is the minimal value for
which $F^*(x,y,z)\in S$.  If there is no such $z$ then 
$F^*_S(x,y)$ is not defined.  This occurs with probability
$(\frac{m-i}{m})^m$  for any
particular $x,y$ so the probability $F^*_S$ is not always
defined is at most $i^2(\frac{m-i}{m})^m=o(1)$.  
\par We generate a random three-place $F^*$ and then
consider $F^*_S$ with $S$ a random set of size $i=n$ or $n+1$.
Conditioning on $F^*_S$ being always defined it then has
the distribution of a random two-place function on $i$ points.
Thus $\Pr[A]$ over $[n],F_n$ is within $o(1)$ of $\Pr[A]$ when
$F_n=F^*_S$ is chosen in this manner.  Thus, as in \S 2, it
suffices to show for any $F^*$ and $A$ that, letting $g(i)$
denote the probability $F^*_S$ satisfies $A$ with $S$ a
uniformly chosen $i$-set, $g(n+1)-g(n)=o(1)$.  Again fix
$F^*$ and $A$ and let $z_x$ be the Boolean value of
$x\in S$ for $1\leq x\leq 2n+1$.  In $A$ replace the
ternary relation $F(a,b)=c$ by
$\wedge_{y<d} \neg z_{F(a,b,y)}$ \& $z_{F(a,b,d)}$.
(For $d=1$ this is simply
True.)  As in \S 2 replace $\exists_x P(x)$ by $\vee_x
(z_x\wedge P^*(x))$ where  $P^*(x)$ has been inductively
defined as the replacement of $P(x)$.  Then the statement
that $F^*_S$ satisfies $A$ becomes a Boolean function of
the $z_1,\ldots,z_m$, as before it is an $AC^0$ function,
and by \S 2 we have $g(n+1)-g(n)=o(1)$.
\medskip

\centerline{$\ast$\qquad\qquad$\ast$\qquad\qquad$\ast$}

\medskip

The following discussion is directed mainly for logicians but may be of
interest for CS-oriented readers as well.

\begin{discussion}
\label{disc1.6}
{\em
\begin{description}
\item[(1)] Note that the results of \cite{Sh463} cannot be gotten in this way
as the proof here use high symmetry.
The problem there was: let $\bar p=\langle p_i: i\in \bbn\rangle$ be a
sequence of probabilities such that $\sum_i p_i < \infty$.
Let $G(n, \bar p)$ be the random graph with set of nodes $[n]=\{1,
\ldots, n\}$ and the edges drawn independently, and for $i\neq j$ the
probability of $\{i, j\}$ being an edge is $p_{|i-j|}$.

The very weak 0-1 law was proved for this context in \cite{Sh463}
(earlier on this context (probability depending on distance)
was introduced and investigated in \L uczak and
Shelah \cite{LuSh435}). Now drawing $G(2n+1, \bar p)$ and then
restricting ourselves to a random $S\subseteq \{1, \ldots, 2n+1\}$ with
$n$, and with $n+1$ elements, fail as $G(2n+1, \bar p)\restriction_S$
does not have the same distribution as $G(|S|, \bar p)$.

\item[(2)] We may want to phrase the result generally; 
\begin{description}
\item[\ ] One way: just say that
$M_n,M_{n+1}$ can be gotten as above : draw a model on $[2n+1]=\{1,
\ldots, 2n+1\}$ (i.e. with
this universe), then choose randomly subsets $P^\ell_n$ with $n +\ell$
elements and restrict yourself to it. 
\end{description}
\item[(3)] Two random linear order satisfies the very weak 0-1 zero law (mean:
take two random functions from $[n]$ to $[0,1]_{\Bbb R}$). The proof
should be clear.
\item[(4)] All this is for fixed probabilities;  we then can allow
probabilities depending on $n$ e.g. we may consider $G_{<}(n, p_n)$ is
the model with set of elements $\{1, \ldots, n\}$, the order relation
and we draw edges with edge probability $p_n$  depending on $n$.
This call for estimating two number (for $\varphi$ first
order sentence):
$$
\alpha_n=|\Prob(G_<(n, p_{n+1})\models \varphi) -\Prob(G_<(n,
p_n)\models \varphi)|
$$
$$
\beta_n=|\Prob(G_< (n, p_{n+1})\models \varphi) - \Prob(G_< (n+1,
p_{n+1}\models \varphi)|
$$
As for $\beta_n$ the question is how much does the proof here depend on
having the probability $\frac{1}{2}$. Direct inspection on the proof
show it does not at all (this just influence on determining the specific
Boolean function with 2n+1 variables) so we know that $\beta_n$ converge
to zero.

As for $\alpha_n$, clearly the question is how fast $p_n$ change.

\item[(5)] As said in \cite{Sh463} we can also consider
$\Lim\big(\Prob_{n+h(n)} (M_{n+h(n)}\models \psi)- \Prob_n(M_n\models
\psi)\big)=0$, i.e. characterize the function $h$ for which this
holds but this was not dealt with there. Hopefully there is a
threshold phenomena. Probably this family 
of problems will appeal to mathematicians with an analytic
background. 

Another problem, closer to my heart,  
is to understand the model theory: in some
sense first order formulas cannot express too much, but can
we find  a more direct statement fulfilling this?

\centerline{$\ast$\qquad\qquad$\ast$\qquad\qquad$\ast$}

Another way to present the first
problem for our case is: close (or at least narrow) the analytic gap
between \cite{CHSh245} and the present paper.

After this work, Boppana and Spencer \cite{BS}, continuing the present
paper and 
\cite{CHSh245}, address the problem and completely solve it. More
specifically they proved the following.

\begin{quotation}
For every sentence~$A$ there exists a number~$t$ so that
$m(n)=O(n\ln^{-t}n)$ implies
\[ \lim\limits_{n\rightarrow\infty} f_A(n+m(n)) - f_A(n) = 0. \]
\end{quotation}
And
\begin{quotation}
For every number~$t$ there exists a sentence~$A$ and a
function~$m(n)=O(n\ln^{-t}n)$ so that $f_A(n+m(n))-f_A(n)$ does not
approach zero.
\end{quotation}
Together we could say:  a function
$m(n)$ has the property that for all $A$ and all $m'(n)\leq m(n)$ we have
$f_A(n+m'(n))-f_A(n) \rightarrow 0$  if and only if $m(n)=o(n\ln^{-t}n)$ for
all $t$.

For improving the bound from this side they have used Hastad
switching lemma [Hastad] (see [AS], \S11.2, Lemma 2.1).

\item[(6)] If we use logic stronger than first order , it cannot be too strong 
(on monadic logic see \cite{KfSh201}), but we may allow quantification over
subsets of size $k_n$, e.g. $\log(n)$ there are two issues:
\begin{description}
\item[(A)] when for both $n$ and $n+1$ we quantify over subsets of size $k_n$,
we should just increase $M$ by having the set $[n]^{k_n}$ as a set of extra
elements, so in (*), $P$ is chosen as a random subset of $\{1,2,3,\ldots,2n,
2n+1\}$ 
with $n$ or $n+1$ elements but the model has about $(2n+1)^{k_n}$ elements;
this require stronger theorem, still true (up to very near to exponentiation)
\item[(B)] if $k_n \not=  k_{n+1}$  we need to show it does
not matter, 
we may choose to round $k_n = \log_2(n)$ so only for rare $n$ the value
change so we weaken a little the theorem {\em or} we may look at
sentences for which 
this does not matter . 

Maybe more naturally, together with choosing randomly $\cM_n$ we choose
a number $\underline{k}_n$, and the probability  of
$\underline{k}_n=k_n+i$ if $i\in
[-k_n/2, k_n/2]$ being $1/k_n$.

And we ask for ``$p^\varphi_n =: \Prob(\cM_n \models
\varphi$ where the
monadic quantifier is interpreted as varying on set with $\leq
\underline{k}_n$ elements) for sentence $\varphi$ (the point of the
distribution of $\underline{k}_n$ is just that for $n$, $n+1$ they
differ a little).
E.g. if for a random graph on $n$ (probability 0.5) we ask on
the property "the size of maximal clique of size at most $[\log_2 n ]^2$ is
even" it satisfies the very weak zero one law
\end{description}
Of course we know much more on this, still it shows that this old result (more
exactly - a weakened version) can be put in our framework.
\end{description}
}
\end{discussion}

\shlhetal  
\end{document}